\documentclass[12pt]{amsart}
\def\Image{\operatorname{Im}}
\def\N{{\Bbb N}}
\def\Z{{\Bbb Z}}

\newtheorem{Theorem}{Theorem}[section]
\newtheorem{Corollary}[Theorem]{Corollary}
\newtheorem{Lemma}[Theorem]{Lemma}

\theoremstyle{definition}

\newtheorem*{Example}{Example}
\newtheorem*{Problem}{Problem}
\theoremstyle{remark}
\newtheorem*{Remark}{Remark}
\newtheorem*{Ack}{Acknowledgments}
\begin{document}
\sloppy
\title{On boundaries of parabolic subgroups of Coxeter groups}
\author{Tetsuya Hosaka} 
\address{Department of Mathematics, Utsunomiya University, 
Utsunomiya, 321-8505, Japan}
\date{October 4, 2004}
\email{hosaka@cc.utsunomiya-u.ac.jp}
\keywords{boundaries of Coxeter groups}
\subjclass[2000]{57M07, 20F65, 20F55}
\thanks{Partly supported by the Grant-in-Aid for Scientific Research, 
The Ministry of Education, Culture, Sports, Science and Technology, Japan, 
(No.~15740029).}
\maketitle
\begin{abstract}
In this paper, 
we investigate boundaries of parabolic subgroups of Coxeter groups.
Let $(W,S)$ be a Coxeter system and 
let $T$ be a subset of $S$ such that the parabolic subgroup $W_T$ is infinite.
Then we show that if a certain set is quasi-dense in $W$, 
then $W \partial\Sigma(W_T,T)$ is dense in the boundary $\partial\Sigma(W,S)$ of 
the Coxeter system $(W,S)$, 
where $\partial\Sigma(W_T,T)$ is the boundary of $(W_T,T)$.
\end{abstract}

\section{Introduction and preliminaries}

The purpose of this paper is to study 
boundaries of parabolic subgroups of Coxeter groups.
A {\it Coxeter group} is a group $W$ having a presentation
$$\langle \,S \, | \, (st)^{m(s,t)}=1 \ \text{for}\ s,t \in S \,
\rangle,$$ 
where $S$ is a finite set and 
$m:S \times S \rightarrow \N \cup \{\infty\}$ is a function 
satisfying the following conditions:
\begin{enumerate}
\item[(1)] $m(s,t)=m(t,s)$ for each $s,t \in S$,
\item[(2)] $m(s,s)=1$ for each $s \in S$, and
\item[(3)] $m(s,t) \ge 2$ for each $s,t \in S$
such that $s\neq t$.
\end{enumerate}
The pair $(W,S)$ is called a {\it Coxeter system}.
Let $(W,S)$ be a Coxeter system.
For a subset $T \subset S$, 
$W_T$ is defined as the subgroup of $W$ generated by $T$, 
and called a {\it parabolic subgroup}.
If $T$ is the empty set, then $W_T$ is the trivial group.
A subset $T\subset S$ is called a {\it spherical subset of $S$}, 
if the parabolic subgroup $W_T$ is finite.
For each $w \in W$, 
we define $S(w)= \{s \in S \,|\, \ell(ws) < \ell(w)\}$, 
where $\ell(w)$ is the minimum length of 
words in $S$ which represents $w$.
For a subset $T \subset S$, 
we also define $W^T= \{w \in W \,|\, S(w)=T \}$. 

A Coxeter system $(W,S)$ is said to be {\it irreducible} 
if for every nonempty and proper subset $T$ of $S$, 
$W$ does not decompose as 
the direct product of $W_T$ and $W_{S \setminus T}$.
Let $(W,S)$ be a Coxeter system. 
Then there exists a unique decomposition $\{S_1,\ldots,S_r\}$ of $S$ 
such that $W$ is the direct product of 
the parabolic subgroups $W_{S_1},\ldots,W_{S_r}$ and 
each Coxeter system $(W_{S_i},S_i)$ is irreducible 
(cf.\ \cite{Bo}, \cite[p.30]{Hu}).
Here we define $\tilde{S}=\bigcup \{S_i\,|\,W_{S_i}\ \text{is infinite}\}$.
It was proved in \cite{Ho} that 
the parabolic subgroup $W_{\tilde{S}}$ is the minimal parabolic subgroup 
of finite index in $W$.
The parabolic subgroup $W_{\tilde{S}}$ is called 
the {\it essential parabolic subgroup} of $(W,S)$.

Every Coxeter system $(W,S)$ determines 
a {\it Davis-Moussong complex} $\Sigma(W,S)$ 
which is a CAT(0) geodesic space (\cite{D1}, \cite{D2}, \cite{D3}, \cite{M}).
Here the $1$-skeleton of $\Sigma(W,S)$ is 
the Cayley graph of $W$ with respect to $S$.
The natural action of $W$ on $\Sigma(W,S)$ is proper, cocompact and by isometry.
We can consider a certain fundamental domain $K(W,S)$ 
which is called a {\it chamber} of $\Sigma(W,S)$ 
such that $W K(W,S)=\Sigma(W,S)$ (\cite{D2}, \cite{D3}).
If $W$ is infinite, then $\Sigma(W,S)$ is noncompact and 
$\Sigma(W,S)$ can be compactified by adding its ideal boundary
$\partial \Sigma(W,S)$ (\cite{BH}, \cite[\S 4]{D2}, \cite{GH}).
This boundary 
$\partial \Sigma(W,S)$ is called the {\it boundary of} $(W,S)$.
We note that the natural action of $W$ on $\Sigma(W,S)$ 
induces an action of $W$ on $\partial \Sigma(W,S)$.
For each subset $T\subset S$, 
$\Sigma(W_T,T)$ is a subcomplex of $\Sigma(W,S)$ and 
the boundary $\partial\Sigma(W_T,T)$ of $(W_T,T)$ is a subspace 
of $\partial \Sigma(W,S)$.

A subset $A$ of a space $X$ is said to be {\it dense} in $X$, 
if $\overline{A}=X$. 
A subset $A$ of a metric space $X$ is said to be {\it quasi-dense}, 
if there exists $N>0$ such that 
each point of $X$ is $N$-close to some point of $A$.

Let $(W,S)$ be a Coxeter system. 
Then $W$ has the {\it word metric} $d_{\ell}$ 
defined by $d_{\ell}(w,w')=\ell(w^{-1}w')$ 
for each $w,w'\in W$. 

We denote by $o(g)$ the order of 
an element $g$ in a Coxeter group $W$.

After some preliminaries in Section~2, 
we prove the following theorem in Section~3.

\begin{Theorem}\label{Thm}
Let $(W,S)$ be a Coxeter system and 
let $T$ be a subset of $S$ such that $W_T$ is infinite.
If the set 
\begin{align*}
\bigcup\{W^{\{s\}}\,|\, s\in S \ \text{such that}&\ 
o(ss_0)=\infty \ \text{and} \ s_0t\neq ts_0 \\
&\ \text{for some}\ s_0\in S\setminus T \ \text{and}\ t\in \tilde{T}\}
\end{align*}
is quasi-dense in $W$ with respect to the word metric, 
then $W \partial\Sigma(W_T,T)$ 
is dense in $\partial\Sigma(W,S)$, 
where $W_{\tilde{T}}$ is 
the essential parabolic subgroup of $(W_T,T)$.
\end{Theorem}

\begin{Remark}
For a Gromov hyperbolic group $G$ and 
the boundary $\partial G$ of $G$, 
we can show that $G\alpha$ is dense in $\partial G$ 
for any $\alpha\in\partial G$ by an easy argument.
Hence 
if $W$ is a hyperbolic Coxeter group, 
then $W\partial\Sigma(W_T,T)$ is dense in $\partial\Sigma(W,S)$ 
for any $T\subset S$ such that $W_T$ is infinite.
\end{Remark}

As an application of Theorem~\ref{Thm}, 
we obtain the following corollary.

\begin{Corollary}\label{Cor}
Let $(W,S)$ be a Coxeter system and
let $T$ be a subset of $S$ such that $W_T$ is infinite.
Suppose that 
there exist a maximal spherical subset $U$ of $S$ 
and an element $s\in S$ such that $o(su)\ge 3$ for every $u\in U$
and $o(su_0)=\infty$ for some $u_0\in U$.
If 
\begin{enumerate}
\item[(1)] $s\not\in T$ and $u_0\in \tilde{T}$, or
\item[(2)] $u_0\not\in T$ and $s\in \tilde{T}$,
\end{enumerate}
then $W \partial\Sigma(W_T,T)$ 
is dense in $\partial\Sigma(W,S)$.
\end{Corollary}

\begin{Example}
We consider the Coxeter system $(W,\{s_1,s_2,s_3,s_4\})$ 
defined by the diagram in Figure~1, 
where Figure~1 defines the Coxeter function $m$ 
for $(W,\{s_1,s_2,s_3,s_4\})$ 
and is not a classical Coxeter diagram.
We note that the Coxeter group $W$ 
is not hyperbolic in the Gromov sense, 
since it contains a copy of $\Z^2$.

\begin{figure}[htbp]
\unitlength = 0.9mm
\begin{center}
\begin{picture}(80,28)(-40,-5)
\put(0,0){\line(1,0){24}}
\put(0,0){\line(-1,0){24}}
\put(0,0){\line(-2,3){12}}
\put(-24,0){\line(2,3){12}}
\put(0,0){\circle*{1.3}}
\put(24,0){\circle*{1.3}}
\put(-24,0){\circle*{1.3}}
\put(-12,18){\circle*{1.3}}
{\small
\put(11,-4){$3$}
\put(-13,-4){$3$}
\put(-5,9){$3$}
\put(-21,9){$3$}
\put(-13.5,20){$s_1$}
\put(-26,-4){$s_2$}
\put(-1.5,-4){$s_3$}
\put(23,-4){$s_4$}
}
\end{picture}
\end{center}
\caption[]{}
\label{fig1}
\end{figure}

Let $T=\{s_1,s_4\}$.
Then we put $U=\{s_3,s_4\}$, $s=s_2$ and $u_0=s_4$.
Since $s\not\in T$ and $u_0\in \tilde{T}$, 
by Corollary~\ref{Cor}~(1), 
$W \partial\Sigma(W_T,T)$ 
is dense in $\partial\Sigma(W,S)$.
By the same argument, 
if $T=\{s_2,s_4\}$ then 
$W \partial\Sigma(W_T,T)$ 
is also dense in $\partial\Sigma(W,S)$.
This implies that 
$W \partial\Sigma(W_T,T)$ 
is dense in $\partial\Sigma(W,S)$
for each $T\in\{\{s_1,s_3,s_4\},\{s_1,s_2,s_4\},\{s_2,s_3,s_4\}\}$.

Let $T=\{s_1,s_2,s_3\}$.
Then we put $U=\{s_3,s_4\}$, $s=s_1$ and $u_0=s_4$.
Since $u_0\not\in T$ and $s\in \tilde{T}$, 
by Corollary~\ref{Cor}~(2), 
$W \partial\Sigma(W_T,T)$ 
is dense in $\partial\Sigma(W,S)$.

Thus, in this example, 
$W \partial\Sigma(W_T,T)$ 
is dense in $\partial\Sigma(W,S)$ 
for any subset $T\subset S$ such that $W_T$ is infinite.
\end{Example}

Concerning $W$-invariantness of $\partial\Sigma(W_T,T)$, 
the following theorem is known.

\begin{Theorem}[\cite{Ho}]~
\begin{enumerate}
\item[(1)] Let $(W,S)$ be a Coxeter system and $T\subset S$. 
Then $\partial \Sigma(W_T,T)$ is $W$-invariant if and only if 
$W=W_{\tilde{T}}\times W_{S\setminus \tilde{T}}$.
\item[(2)] Let $(W,S)$ be an irreducible Coxeter system and
let $T$ be a proper subset of $S$ such that $W_T$ is infinite.
Then $\partial \Sigma(W_T,T)$ is not $W$-invariant.
\end{enumerate}
\end{Theorem}

The following problem is open.

\begin{Problem}
Let $(W,S)$ be a Coxeter system and
let $T$ be a subset of $S$ such that $W_T$ is infinite.
Is it the case that 
if $\partial\Sigma(W_T,T)$ is not $W$-invariant then
$W\partial\Sigma(W_T,T)$ is dense in $\partial\Sigma(W,S)$? 
Particulary, 
is it the case that 
if $(W,S)$ is an irreducible Coxeter system then
$W\partial\Sigma(W_T,T)$ is dense in $\partial\Sigma(W,S)$ 
for any subset $T$ of $S$ such that $W_T$ is infinite? 
\end{Problem}

\begin{Ack}
The author would like to thank the referee for helpful advice.
\end{Ack}

\section{Lemmas on Coxeter groups and $\Sigma(W,S)$}

In this section, we recall and prove some lemmas for Coxeter groups 
and $\Sigma(W,S)$ which are used later.

Let $(W,S)$ be a Coxeter system.
For a subset $T$ of $S$, we define 
$A_T=\{w \in W\,|\,\ell(wt)>\ell(w)\ \text{for all}\ t\in T\}$.
For $w\in W$, 
a representation $w=s_1\cdots s_l$ ($s_i \in S$) is said to be 
{\it reduced}, if $\ell(w)=l$.

The following lemmas are known.

\begin{Lemma}[{\cite[Lemma~2.4]{Ho}}]\label{lem1-2.4}
Let $(W,S)$ be a Coxeter system and $T \subset S$.
Then $[W:W_T]=|A_T|$.
\end{Lemma}

\begin{Lemma}[{\cite[Lemma~2.7]{Ho}}]\label{lem1-2.7}
Let $(W,S)$ be a Coxeter system, 
$T \subset S$ and $t_1,\ldots,t_n \in S\setminus T$ 
($t_i\neq t_j$ if $i\neq j$).
Suppose that 
$t_i t_{i+1} \neq t_{i+1}t_i$
for any $1\le i \le n-1$ and 
$t_n t \neq t t_n$
for any $t \in T$.
Then 
$$(W_{S\setminus \{t_1,\ldots, t_n\}})^T t_n\cdots t_1
\subset W^{\{t_1\}}.$$
\end{Lemma}

\begin{Theorem}[{\cite[Corollary~3.4]{Ho}}]\label{thm1-3.4}
Let $(W,S)$ be a Coxeter system and $T \subset S$.
If $W_T$ is a parabolic subgroup of finite index in $W$, then
\begin{enumerate}
\item[$(1)$] $\tilde{S} \subset T$,
\item[$(2)$] $W_T=W_{\tilde{S}} \times W_{T \setminus \tilde{S}}$ and
\item[$(3)$] $[W:W_T]=|W_{S\setminus \tilde{S}}|/|W_{T\setminus \tilde{S}}|$.
\end{enumerate}
\end{Theorem}

\begin{Lemma}[{\cite[Lemma~2.4]{Ho2}}]\label{lem3}
Let $(W,S)$ be a Coxeter system, $w\in W$ and $s_0\in S$. 
Suppose that $o(s_0t)\ge 3$ for every $t\in S(w)$ and 
that $o(s_0t_0)=\infty$ for some $t_0\in S(w)$.
Then $ws_0\in W^{\{s_0\}}$.
\end{Lemma}

\begin{Lemma}[{\cite[Lemma~2.5]{Ho2}}]\label{lem4}
Let $(W,S)$ be a Coxeter system. 
Suppose that there exist a maximal spherical subset $T$ of $S$ 
and $s_0\in S$ such that $o(s_0t)\ge 3$ for every $t\in T$
and $o(s_0t_0)=\infty$ for some $t_0\in T$. 
Then $W^{\{s_0\}}$ is quasi-dense in $W$.
\end{Lemma}

Using Lemmas~\ref{lem1-2.4} and \ref{lem1-2.7}
and Theorem~\ref{thm1-3.4}, 
we prove the following lemma.

\begin{Lemma}\label{Lem1}
Let $(W,S)$ be a Coxeter system, 
let $T$ be a proper subset of $S$ such that $W_T$ is infinite, and 
let 
$$U=\{s\in S\setminus T\,|\, 
W^{\{s\}}s\cap W_T \ \text{is finite}\}.$$
Then $W_{\tilde{T}\cup U}=W_{\tilde{T}}\times W_U$.
\end{Lemma}

\begin{proof}
We note that $S(w)\subset T$ for $w\in W_T$.
Let $u_0\in U$ and 
let $T(u_0)=\{t\in T\,|\, tu_0\neq u_0t\}$.
We first show that 
$W_{T\setminus T(u_0)}$ is a subgroup of finite index in $W_T$.
Here we note that 
$[W_T:W_{T\setminus T(u_0)}]=|A_{T\setminus T(u_0)}\cap W_T|$
by Lemma~\ref{lem1-2.4}.
Then
\begin{align*}
\bigcup_{T'\subset T(u_0)}(W_T)^{T'}
&=\{w\in W_T\,|\,S(w)\subset T(u_0)\} \\
&=\{w\in W_T\,|\,T\setminus T(u_0)\subset T\setminus S(w)\} \\
&=A_{T\setminus T(u_0)}\cap W_T.
\end{align*}
We show that $(W_T)^{T'}$ is finite for any $T'\subset T(u_0)$.
Let $T'\subset T(u_0)$.
Since $tu_0\neq u_0t$ for any $t\in T'$, 
$(W_T)^{T'}u_0\subset W^{\{u_0\}}$ by Lemma~\ref{lem1-2.7}.
Hence $(W_T)^{T'}\subset W^{\{u_0\}}u_0\cap W_T$ which is finite 
because $u_0\in U$.
Thus $(W_T)^{T'}$ is finite for any $T'\subset T(u_0)$, 
and $[W_T:W_{T\setminus T(u_0)}]=|A_{T\setminus T(u_0)}\cap W_T|$ 
is finite.
By Theorem~\ref{thm1-3.4}, 
$\tilde{T}\subset T\setminus T(u_0)$.
Hence $T(u_0)\subset T\setminus \tilde{T}$ for any $u_0\in U$.
Let $A=\{t\in T\,|\,tu_0\neq u_0t \ \text{for some}\ u_0\in U\}$.
Then 
$A=\bigcup_{u_0\in U}T(u_0) \subset T\setminus \tilde{T}$ and 
$$\tilde{T}\subset T\setminus A=
\{t\in T\,|\, tu=ut\ \text{for every}\ u\in U\}.$$
Thus $tu=ut$ for any $t\in\tilde{T}$ and $u\in U$.
This means that $W_{\tilde{T}\cup U}=W_{\tilde{T}}\times W_U$.
\end{proof}

We immediately obtain the following lemma from Lemma~\ref{Lem1}.

\begin{Lemma}\label{Lem2}
Let $(W,S)$ be a Coxeter system, 
let $T$ be a proper subset of $S$ such that $W_T$ is infinite and 
let $s\in S\setminus T$.
If $st\neq ts$ for some $t\in\tilde{T}$, 
then $W^{\{s\}}s\cap W_T$ is infinite.
\end{Lemma}

Concerning $\Sigma(W,S)$, the following lemma is known.

\begin{Lemma}[{\cite[Lemma~3.3]{Ho2}}]\label{lem6}
Let $(W,S)$ be a Coxeter system and $x,y\in W$.
If $o(st)=\infty$ for each $s\in S(x)$ and $t\in S(y^{-1})$, 
then $d(x,\Image \xi_{xy}) \le N$, 
where $\xi_{xy}$ is the geodesic from $1$ to $xy$ in $\Sigma(W,S)$
and $N$ is the diameter of $K(W,S)$ in $\Sigma(W,S)$.
\end{Lemma}

\section{Proof of the main results}

Using the lemmas of Section~2, 
we prove the main results.

\begin{proof}[Proof of Theorem~\ref{Thm}]
Suppose that 
\begin{align*}
A:=\bigcup\{W^{\{s\}}\,|\, s\in S \ \text{such that}&\ 
o(ss_0)=\infty \ \text{and} \ s_0t\neq ts_0 \\
&\ \text{for some}\ s_0\in S\setminus T \ \text{and}\ t\in \tilde{T}\}
\end{align*}
is quasi-dense in $W$. 

We first show that 
for each $w\in A$, 
there exists $v\in W$ and $\alpha\in\partial\Sigma(W_T,T)$ such that 
$d(w,\Image\xi_{v\alpha})\le N$, 
where $N$ is the diameter of $K(W,S)$ in $\Sigma(W,S)$ 
and $\xi_{v\alpha}$ is the geodesic ray issuing from $1$ 
such that $\xi_{v\alpha}(\infty)=v\alpha$. 

Let $w\in A$. 
Then $w\in W^{\{s\}}$, $o(ss_0)=\infty$ and $s_0t\neq ts_0$
for some $s\in S$, $s_0\in S\setminus T$ and $t\in\tilde{T}$. 
By Lemma~\ref{Lem2}, 
$W^{\{s_0\}}s_0\cap W_T$ is infinite.
Hence 
there exists a sequence 
$\{x_i\}\subset (W^{\{s_0\}}s_0\cap W_T)^{-1}$
which converges to some point 
$\alpha\in\partial\Sigma(W_T,T)$.
Since $x_i\in(W^{\{s_0\}}s_0)^{-1}$, 
$(s_0x_i)^{-1}=x_i^{-1}s_0\in W^{\{s_0\}}$.
By Lemma~\ref{lem6}, 
$d(w,\Image\xi_{ws_0x_i})\le N$ for any $i$ 
because $w\in W^{\{s\}}$, $s_0x_i\in(W^{\{s_0\}})^{-1}$ 
and $o(ss_0)=\infty$.
Hence $d(w,\Image\xi_{ws_0\alpha})\le N$.

For each $\beta\in \partial\Sigma(W,S)$, 
there exists a sequence 
$\{w_i\}\subset A$ 
which converges to $\beta$, 
because $A$ is quasi-dense in $W$. 
By the above argument, 
there exist sequences $\{v_i\}\subset W$ 
and $\{\alpha_i\}\subset\partial\Sigma(W_T,T)$ such that 
$d(w_i,\Image\xi_{v_i\alpha_i})\le N$ for each $i$.
Then the sequence $\{v_i\alpha_i\}$ converges to $\beta$ in 
$\partial\Sigma(W,S)$ 
because $\{w_i\}$ converges to $\beta$. 
Therefore 
$W \partial\Sigma(W_T,T)$ is dense in $\partial\Sigma(W,S)$.
\end{proof}

\begin{proof}[Proof of Corollary~\ref{Cor}]
Suppose that 
there exist a maximal spherical subset $U$ of $S$ 
and an element $s\in S$ such that $o(su)\ge 3$ for any $u\in U$
and $o(su_0)=\infty$ for some $u_0\in U$.
Then $W^{\{s\}}$ is quasi-dense in $W$ by Lemma~\ref{lem4}.

(1) If $s\not\in T$ and $u_0\in \tilde{T}$, 
then $W^{\{u_0\}}$ is quasi-dense in $W$ 
because $W^{\{s\}}u_0\subset W^{\{u_0\}}$ by Lemma~\ref{lem3}.
Hence $W \partial\Sigma(W_T,T)$ 
is dense in $\partial\Sigma(W,S)$ by Theorem~\ref{Thm}.

(2) If $u_0\not\in T$ and $s\in \tilde{T}$, 
then by Theorem~\ref{Thm}, 
$W \partial\Sigma(W_T,T)$ is dense in $\partial\Sigma(W,S)$, 
because $o(su_0)=\infty$, $u_0\in S\setminus T$ and $s\in \tilde{T}$.
\end{proof}

%

%
\end{document}